\RequirePackage{etoolbox}
\csdef{input@path}{{style/}{graphics/}}
\documentclass[aos,MSNbibl,dvips]{arximspdf}
\usepackage{mathbh}

\doi{10.1214/14-AOS1293}
\volume{43}
\issue{2}
\pubyear{2015}
\firstpage{592}
\lastpage{619}
\docsubty{FLA}

\makeatletter
\newcommand{\rrvert}{\vert}
\newcommand{\rrVert}{\Vert}
\newcommand{\llvert}{\vert}
\newcommand{\llVert}{\Vert}
\renewcommand{\mid}{|}
\newcommand{\Vo}{\mathop{V}^{\circ}\hspace*{-3pt}{}}
\newcommand{\underset}[2]{\mathop{#2}\limits_{#1}}
\newtheorem{teo}{Theorem}[section]
\newtheorem{lem}{Lemma}[section]
\newproclaim{Def}{Definition}[section]
\newproclaim{ex}{Example}
\newtheorem{prop}{Proposition}[section]
\renewcommand{\c}{\mathbf{c}}
\def\Var{\operatorname{Var}}
\makeatother

\begin{document}
\begin{frontmatter}

\title{Nonasymptotic bounds for vector quantization in~Hilbert spaces}
\runtitle{Nonasymptotic bounds for vector quantization}

\begin{aug}
\author[A]{\fnms{Cl\'{e}ment}~\snm{Levrard}\corref{}\ead[label=e1]{clement.levrard@inria.fr}}
\runauthor{C. Levrard}
\affiliation{Universit\'{e} Paris Sud, UPMC and INRIA}
\address[A]{B\^atiment Alan Turing\\
Campus de l\'Ecole Polytechnique\\
INRIA\\
91120 Palaiseau\\
France\\
\printead{e1}}
\end{aug}

%
\received{\smonth{11} \syear{2013}}
%
\revised{\smonth{11} \syear{2014}}

%
\begin{abstract}
Recent results in quantization theory show that the mean-squared
expected distortion can reach a rate of convergence of $\mathcal{O}(1/n)$, where $n$ is
the sample size [see, e.g.,
\textit{IEEE Trans. Inform. Theory} \textbf{60} (2014) 7279--7292
or
\textit{Electron. J. Stat.} \textbf{7} (2013) 1716--1746].
This rate is attained
for the empirical risk minimizer strategy, if the source distribution
satisfies some
regularity conditions. However, the dependency of the average
distortion on other
parameters is not known, and these results are only valid for
distributions over
finite-dimensional Euclidean spaces.

This paper deals with the general case of distributions over separable, possibly
infinite dimensional, Hilbert spaces. A condition is proposed, which
may be thought
of as a margin condition [see, e.g.,
\textit{Ann. Statist.} \textbf{27} (1999) 1808--1829], under which a
nonasymptotic
upper bound on the expected distortion rate of the empirically optimal
quantizer is
derived. The dependency of the distortion on other parameters of
distributions is
then discussed, in particular through a minimax lower bound.
\end{abstract}

%
\end{frontmatter}

\section{Introduction}\label{Introduction}

Quantization, also called lossy data compression in information theory,
is the problem of replacing a probability distribution with an
efficient and compact representation, that is a finite set of points.
To be more precise, let $\mathcal{H}$ denote a separable Hilbert
space, and let $P$ denote a probability distribution over $\mathcal
{H}$. For a positive integer $k$, a so-called $k$-points quantizer $Q$
is a map from $\mathcal{H}$ to $\mathcal{H}$, whose image set is made
of exactly $k$ points, that is $\llvert  Q(\mathcal{H}) \rrvert   = k$.
For such a quantizer, every image point $c_i \in Q (\mathcal
{H}  )$ is called a code point, and the vector composed of the
code points $(c_1, \ldots, c_k)$ is called a codebook, denoted by $\c
$. By considering the pre-images of its code points, a quantizer $Q$
partitions the separable Hilbert space $\mathcal{H}$ into $k$ groups,
and assigns each group a representative. General references on the
subject are to be found in \cite{GL00,Gersho91} and \cite
{Linder02} among others.

The quantization theory was originally developed as a way to answer
signal compression issues in the late 1940s (see, e.g., \cite
{Gersho91}). However, unsupervised classification is also in the scope
of its application. Isolating meaningful groups from a cloud of data is
a topic of interest in many fields, from social science to biology.
Classifying points into dissimilar groups of similar items is more
interesting as the amount of accessible data is large. In many cases
data need to be preprocessed through a quantization algorithm in order
to be exploited.

If the distribution $P$ has a finite second moment, the performance of
a quantizer $Q$ is measured by the risk, or distortion
\[
R(Q):= P \bigl\llVert x - Q(x) \bigr\rrVert ^2,
\]
where $Pf$ means integration of the function $f$ with respect to $P$.
The choice of the squared norm is convenient, since it takes advantages
of the Hilbert space structure of $\mathcal{H}$. Nevertheless, it is
worth pointing out that several authors deal with more general
distortion functions. For further information on this topic, the
interested reader is referred to \cite{GL00} or \cite{Fischer10}.

In order to minimize the distortion introduced above, it is clear that
only quantizers of the type $x \mapsto\arg\min_{c_1, \ldots, c_k}{\llVert
x - c_i\rrVert  ^2}$ are to be considered. Such quantizers are called
nearest-neighbor quantizers. With a slight abuse of notation, $R(\c)$
will denote the risk of the nearest-neighbor quantizer associated with
a codebook $\c$.

Provided that $P$ has a bounded support, there exist optimal codebooks
minimizing the risk $R$ (see, e.g., Corollary 3.1 in \cite{Fischer10}
or Theorem~1 in \cite{Graf07}). The aim is to design a codebook $\hat
{\c}_n$, according to an $n$-sample drawn from $P$, whose distortion
is as close as possible to the optimal distortion $R(\c^*)$, where $\c
^*$ denotes an optimal codebook.

To solve this problem, most approaches to date attempt to implement the
principle of empirical risk minimization in the vector quantization
context. Let $X_1, \ldots, X_n$ denote an independent and identically
distributed sample with distribution $P$. According to this principle,
good code points can be found by searching for ones that minimize the
empirical distortion over the training data, defined by
\[
\hat{R}_n(\mathbf{c}):= \frac{1}{n}\sum
_{i=1}^{n}{\bigl\llVert X_i -
Q(X_i)\bigr\rrVert ^2} = \frac{1}{n}\sum
_{i=1}^{n}{ \min_{j=1, \ldots,
k}{\llVert
X_i - c_j\rrVert ^2}}.
\]
If the training data represents the source well, then $\hat{\mathbf
{c}}_n$ will hopefully also perform near optimally on the real source,
that is, $\ell(\hat{\c}_n,\c^*) = R(\hat{\c}_n)-R(\c^*) \approx
0$. The problem of quantifying how good empirically designed codebooks
are, compared to the truly optimal ones, has been extensively studied,
as, for instance, in \cite{Linder02} in the finite-dimensional case.

If $\mathcal{H} = \mathbb{R}^d$, for some $d>0$, it has been proved
in \cite{Linder94} that $\mathbb{E}\ell(\hat{\c}_n,\c^*) =
\mathcal{O}(1/\sqrt{n})$, provided that $P$ has a bounded support.
This result has been extended to the case where $\mathcal{H}$ is a
separable Hilbert space in \cite{Biau08}. However, this upper bound
has been tightened whenever the source distribution satisfies
additional assumptions, in the finite-dimensional case only.

When $\mathcal{H} = \mathbb{R}^d$, for the special case of finitely
supported distributions, it is shown in \cite{Antos04} that $\mathbb
{E}\ell(\hat{\c}_n,\c^*) = \mathcal{O}(1/n)$. There are much more
results in the case where $P$ is not assumed to have a finite support.

In fact, different sets of assumptions have been introduced in \cite
{Antos04,Pollard82} or \cite{Levrard12}, to derive fast
convergence rates for the distortion in the finite-dimensional case. To
be more precise, it is proved in \cite{Antos04} that, if $P$ has a
support bounded by $M$ and satisfies a technical inequality, namely for
some fixed $a>0$, for every codebook~$\c$, there is a $\c^*$ optimal
codebook such that
%
\begin{equation}
\label{Antoscondition} \ell\bigl(\c,\c^*\bigr) \geq a \Var \Bigl( \min
_{j=1, \ldots, k}{\llVert X - c_j\rrVert ^2} -
\min_{j=1, \ldots, k}{\bigl\llVert X - c^*_j\bigr\rrVert
^2} \Bigr),
\end{equation}
then $\mathbb{E} \ell(\hat{\c}_n,\c^*) \leq C(k,d,P) \log(n)/n$,
where $C(k,d,P)$ depends on the natural parameters $k$ and $d$, and
also on $P$, but only through $M$ and the technical parameter $a$.
However, in the continuous density and unique minimum case, it has been
proved in \cite{Chou94}, following the approach of \cite{Pollard82},
that provided the Hessian matrix of $\c\mapsto R(\c)$ is positive
definite at the optimal codebook, $n \ell(\hat{\c}_n,\c^*)$
converges in distribution to a law, depending on the Hessian matrix. As
proved in \cite{Levrard12}, the technique used in \cite{Pollard82}
can be slightly modified to derive a nonasymptotic bound of the type
$\mathbb{E} \ell(\hat{\c}_n,\c^*) \leq C/n$ in this case, for some
unknown $C>0$.

As shown in \cite{Levrard12}, these different sets of assumptions turn
out to be equivalent in the continuous density case to a technical
condition, similar to that used in \cite{Massart06} to derive fast
rates of convergence in the statistical learning framework.

Thus, a question of interest is to know whether some margin type
conditions can be derived for the source distribution to satisfy the
technical condition mentioned above, as has been done in the
statistical learning framework in \cite{Tsybakov99}. This paper
provides a condition, which can clearly be thought of as a margin
condition in the quantization framework, under which condition (\ref
{Antoscondition}) is satisfied. The technical constant $a$ has then an
explicit expression in terms of natural parameters of $P$ from the
quantization point of view. This margin condition does not require
$\mathcal{H}$ to have a finite dimension, or $P$ to have a continuous
density. In the finite-dimensional case, this condition does not demand
either that there exists a unique optimal codebook, as required in
\cite{Pollard82}, hence seems easier to check.

Moreover, a nonasymptotic bound of the type $\mathbb{E} \ell(\hat
{\c}_n,\c^*) \leq C(k,P)/n$ is derived for distributions satisfying
this margin condition, where $C(k,P)$ is explicitly given in terms of
parameters of $P$. This bound is also valid if $\mathcal{H}$ is
infinite dimensional. This point may be of interest for curve
quantization, as done in \cite{Fischer12}.

In addition, a minimax lower bound is given which allows one to discuss
the influence of the different parameters mentioned in the upper bound.
It is worth pointing out that this lower bound is valid over a set of
probability distributions with uniformly bounded continuous densities
and unique optimal codebooks, such that the minimum eigenvalues of the
second derivative matrices of the distortion, at the optimal codebooks,
are uniformly lower bounded. This result generalizes the previous
minimax bound obtained in Theorem~4 of \cite{Antos05} for $k\geq3$
and $d>1$.

This paper is organized as follows. In Section~\ref{Notation}, some
notation and definitions are introduced, along with some basic results
for quantization in a Hilbert space. The so-called margin condition is
then introduced, and the main results are exposed in Section~\ref
{Results}: first an oracle inequality on the loss is stated, along with
a minimax result. Then it is shown that Gaussian mixtures are in the
scope of the margin condition. Finally, the main results are proved in
Section~\ref{Proofs} and the proofs of several supporting lemmas are
deferred to the supplementary material \cite{supple}.

\section{Notation and definitions}\label{Notation}

Throughout this paper, for $M >0$ and $a$ in $\mathcal{H}$, $\mathcal
{B}(a,M)$ and $\mathcal{B}^o(a,M)$ will denote, respectively, the closed
and open ball with center $a$ and radius $M$. For a subset $A$ of
$\mathcal{H}$, $\bigcup_{a \in A}{\mathcal{B}(a,M)}$ will be denoted
by $\mathcal{B}(A,M)$. With a slight abuse of notation, $P$ is said to
be $M$-bounded if its support is included in $\mathcal{B}(0,M)$.
Furthermore, it will also be assumed that the support of $P$ contains
more than $k$ points.

To frame quantization as an empirical risk minimization issue, the
following contrast function $\gamma$ is introduced as
\begin{eqnarray*}
&&\gamma\dvtx  \cases{\displaystyle (\mathcal{H} )^k \times\mathcal{H}
\longrightarrow \mathbb{R},
\vspace*{3pt}\cr
\displaystyle (\mathbf{c},x) \longmapsto \min
_{j=1, \ldots, k}{\llVert x-c_j \rrVert ^2},}
\end{eqnarray*}
where $\c= (c_1, \ldots, c_k)$ denotes a codebook, that is a $k
d$-dimensional vector if \mbox{$\mathcal{H} = \mathbb{R}^d$}. In this paper,
only the case $k \geq2$ will be considered. The risk $R(\c)$ then
takes the form $R(\c)= R(Q) = P\gamma(\c,\cdot)$, where we recall that
$Pf$ denotes the integration of the function $f$ with respect to $P$.
Similarly, the empirical risk $\hat{R}_n(\c)$ can be defined as $\hat
{R}_n(\c) = P_n\gamma(\c,\cdot)$, where $P_n$ is the empirical
distribution associated with $X_1, \ldots, X_n$, in other words
$P_n(A) = (1/n) \llvert   \{i \mid X_i \in A \} \rrvert  $, for any measurable
subset $A \subset\mathcal{H}$.

It is worth pointing out that, if $P$ is $M$-bounded, for some $M >0$,
then there exist such minimizers $\hat{\c}_n$ and $\mathbf{c}^*$
(see, e.g., Corollary 3.1 in \cite{Fischer10}). In the sequel, the set
of minimizers of the risk $R$ will be denoted by $\mathcal{M}$. Since
every permutation of the labels of an optimal codebook provides an
optimal codebook, $\mathcal{M}$ contains more than $k!$ elements. To
address the issue of a large number of optimal codebooks, $\bar
{\mathcal{M}}$ is introduced as a set of codebooks which satisfies
\begin{eqnarray*}
&&\cases{\displaystyle \forall\c^* \in\mathcal{M}, \exists\bar{\c} \in\bar {\mathcal{M}}, &\quad $\bigl\{c_1^*, \ldots, c_k^*\bigr\} = \{
\bar{c}_1, \ldots, \bar{c}_k \}$,
\vspace*{3pt}\cr
\displaystyle \forall\bar{
\c}^1 \neq\bar{\c}^2 \in\bar{\mathcal{M}}, &\quad $\bigl\{\bar{c}^1_1, \ldots, \bar{c}_k^1
\bigr\} \neq \bigl\{ \bar {c}^2_1, \ldots,
\bar{c}_k^2 \bigr\}$.}
\end{eqnarray*}
In other words, $\bar{\mathcal{M}}$ is a subset of the set of optimal
codebooks which contains every element of $\mathcal{M}$, up to a
permutation of the labels, and in which two different codebooks have
different sets of code points. It may be noticed that $\bar{\mathcal
{M}}$ is not uniquely defined. However, when $\mathcal{M}$ is finite,
all the possible $\bar{\mathcal{M}}$ have the same cardinality.

Let $c_1, \ldots, c_k$ be a sequence of code points. A central role is
played by the set of points which are closer to $c_i$ than to any other
$c_j$'s. To be more precise, the Voronoi cell, or quantization cell
associated with $c_i$ is the closed set defined by
\begin{eqnarray*}
\label{Voronoidefinition} V_i(\c)&=& \bigl\{ x \in\mathcal{H} \mid \forall j
\neq i,  \llVert x-c_i\rrVert \leq\llVert x-c_j
\rrVert \bigr\}.
\end{eqnarray*}
Note that $(V_1(\c), \ldots, V_k(\c))$ does not form a partition of
$\mathcal{H}$, since $V_i(\c) \cap V_j(\c) $ may be nonempty. To
address this issue, a Voronoi partition associated with $\c$ is
defined as a sequence of subsets $(W_1(\c), \ldots, W_k(\c))$ which
forms a partition of $\mathcal{H}$, and such that for every $i=1,
\ldots, k$,
\[
\bar{W}_i(\c)=V_i(\c),
\]
where $\bar{W}_i(\c)$ denotes the closure of the subset $W_i(\c)$.
The open Voronoi cell is defined the same way by
\begin{eqnarray*}
\Vo_i(\c)&=& \bigl\{ x \in\mathcal{H} \mid \forall j \neq i, \llVert x-c_i\rrVert < \llVert x-c_j\rrVert \bigr
\}.
\end{eqnarray*}
Given a Voronoi partition $W(\c)=(W_1(\c), \ldots, W_k(\c))$, the
following inclusion holds, for $i$ in $ \{1,\ldots,k \}$,
\[
\Vo_i(\c) \subset W_i(\c) \subset V_i(
\c),
\]
and the risk $R(\c)$ takes the form
\[
R(\c)= \sum_{i=1}^{k}{P \bigl( \llVert
x-c_i\rrVert ^2 \mathbh{1}_{W_i(\c)}(x) \bigr)},
\]
where $\mathbh{1}_A$ denotes the indicator function associated with
$A$. In the case where $(W_1, \ldots, W_k)$ are fixed subsets such
that $P(W_i) \neq0$, for every $i=1, \ldots, k$, it is clear that
\[
P\bigl(\llVert x - c_i\rrVert ^2 \mathbh{1}_{W_i(\c)}(x)
\bigr) \geq P\bigl(\llVert x - \eta_i\rrVert ^2
\mathbh{1}_{W_i(\c)}(x)\bigr),
\]
with equality only if $c_i = \eta_i$, where $\eta_i$ denotes the
conditional expectation of $P$ over the subset $W_i(\c)$, that is,
\[
\eta_i = \frac{P(x\mathbh{1}_{W_i(\c)}(x))}{P(W_i(\c))}.
\]
Moreover, it is proved in Proposition 1 of \cite{Graf07} that, for
every Voronoi partition $W(\c^*)$ associated with an optimal codebook
$\c^*$, and every $i=1, \ldots, k$, $P(W_i(\c^*)) \neq0$.
Consequently, any optimal codebook satisfies the so-called centroid
condition (see, e.g., Section~6.2 of \cite{Gersho91}), that is,
\[
\c^*_i= \frac{P(x\mathbh{1}_{W_i(\c^*)}(x))}{P(W_i(\c^*))}.
\]
As a remark, the centroid condition ensures that $\mathcal{M} \subset
\mathcal{B}(0,M)^k$, and, for every $\c^*$ in $\mathcal{M}$, $i \neq j$,
\begin{eqnarray*}
P\bigl(V_i\bigl(\c^*\bigr) \cap V_j\bigl(\c^*\bigr)
\bigr) & =& P \bigl( \bigl\{ x \in\mathcal{H}\mid \forall i',
\bigl\llVert x-c_i^*\bigr\rrVert = \bigl\llVert x -
c_j^*\bigr\rrVert \leq\bigl\llVert x - c^*_{i'}\bigr
\rrVert \bigr\} \bigr)
\\
&=& 0.
\end{eqnarray*}
A proof of this statement can be found in Proposition 1 of \cite
{Graf07}. According to this remark, it is clear that, for every optimal
Voronoi partition $(W_1(\c^*), \ldots, W_k(\c^*))$,
%
\begin{equation}
\label{nomansland} \cases{\displaystyle P\bigl(W_i\bigl(\c^*\bigr)\bigr) = P
\bigl(V_i\bigl(\c^*\bigr)\bigr),
\vspace*{3pt}\cr
\displaystyle P_n
\bigl(W_i\bigl(\c^*\bigr)\bigr) \underset{\mathrm{a.s.}} {=} P_n
\bigl(V_i\bigl(\c^*\bigr)\bigr).}
\end{equation}
The following quantities are of importance in the bounds exposed in
Section~\ref{mainresult}:
\begin{eqnarray*}
&&\cases{\displaystyle B = \inf_{\c^* \in\mathcal{M}, i \neq j}{\bigl\llVert c_i^* -
c_j^*\bigr\rrVert },
\vspace*{3pt}\cr
\displaystyle p_{\min} = \inf_{\c^* \in\mathcal{M}, i =1,\ldots,k}{P
\bigl(V_i\bigl(\c^*\bigr)\bigr)}.}
\end{eqnarray*}
It is worth noting here that $B \leq2M$ whenever $P$ is $M$-bounded,
and $p_{\min} \leq1/k$. If $\mathcal{M}$ is finite, it is clear that
$p_{\min}$ and $B$ are positive. The following proposition ensures that
this statement remains true when $\mathcal{M}$ is not assumed to be finite.

%
\begin{prop}\label{minorationBetpmin}
Suppose that $P$ is $M$-bounded. Then both $B$ and $p_{\min}$ are positive.
\end{prop}

A proof of Proposition~\ref{minorationBetpmin} is given in
Section~\ref{Proofs}. The role of the boundaries between optimal
Voronoi cells may be compared to the role played by the critical value
$1/2$ for the regression function in the statistical learning framework
(for a comprehensive explanation of this statistical learning point of
view, see, e.g., \cite{Massart06}). To draw this comparison, the
following set is introduced, for any $\c^* \in \mathcal{M}$,
\[
N_{\c^*} = \bigcup_{i \neq j}{V_i
\bigl(\c^*\bigr) \cap V_j\bigl(\c^*\bigr)}.
\]
The region is of importance when considering the conditions under which
the empirical risk minimization strategy for quantization achieves
faster rates of convergence, as exposed in \cite{Levrard12}. However,
to completely translate the margin conditions given in \cite
{Tsybakov99} to the quantization framework, the neighborhood of this
region has to be introduced. For this purpose, the $t$-neighborhood of
the region $N_{\c^*}$ is defined by $\mathcal{B}(N_{\c^*},t)$. The
quantity of interest is the maximal weight of these $t$-neighborhoods
over the set of optimal codebooks, defined by
\[
p(t) = \sup_{\c^* \in\mathcal{M}}P\bigl(\mathcal{B}(N_{\c^*},t)\bigr).
\]
It is straightforward that $p(0)=0$. Intuitively, if $p(t)$ is small
enough, then the source distribution $P$ is concentrated around its
optimal codebook, and may be thought of as a slight modification of the
probability distribution with finite support made of an optimal
codebook $\c^*$. To be more precise, let us introduce the following
key assumption.

\begin{Def}[(Margin condition)]\label{margincondition}
A distribution $P$ satisfies a margin condition with radius $r_0 >0$ if
and only if:
\begin{longlist}[(ii)]
\item[(i)] $P$ is $M$-bounded,
\item[(ii)] for all $0 \leq t \leq r_0$,
%
\begin{equation}
\label{majorationkappa} p(t) \leq\frac{B p_{\min}}{128 M^2}t.
\end{equation}
\end{longlist}
\end{Def}

Note that, since $p(2M) =1$, $p_{\min} \leq1/k$, $k \geq2$ and $B \leq
2M$, (\ref{majorationkappa}) implies that $r_0 < 2M$. It is worth
pointing out that Definition~\ref{margincondition} does not require
$P$ to have a density or a unique optimal codebook, up to relabeling,
contrary to the conditions introduced in \cite{Pollard82}.

Moreover, the margin condition introduced here only requires a local
control of the weight function $p(t)$. The parameter $r_0$ may be
thought of as a gap size around every $N_{\c^*}$, as illustrated by
the following example:

\begin{ex}\label{ex1}
Assume that there exists $r>0$ such that $p(x)=0$
if $x\leq r$ (e.g., if $P$ is supported on $k$ points). Then $P$
satisfies (\ref{majorationkappa}), with radius $r$.
\end{ex}

Note also that the condition mentioned in \cite{Tsybakov99} requires a
control of the weight of the neighborhood of the critical value $1/2$
with a polynomial function with degree larger than $1$. In the
quantization framework, the special role played by the exponent $1$
leads to only consider linear controls of the weight function. This
point is explained by the following example:

\begin{ex}\label{ex2}
Assume that $P$ is $M$-bounded, and that there
exists $Q>0$ and $q>1$ such that $p(x) \leq Q x^q$. Then $P$ satisfies
(\ref{majorationkappa}), with
\[
r_0 = \biggl( \frac{p_{\min} B}{128 M^2 Q} \biggr)^{1/(q-1)}.
\]
\end{ex}

In the case where $P$ has a density and $\mathcal{H} = \mathbb{R}^d$,
the condition (\ref{majorationkappa}) may be considered as a
generalization of the condition stated in Theorem~3.2 of \cite
{Levrard12}, which requires the density of the distribution to be small
enough over every $N_{\c^*}$. In fact, provided that $P$ has a
continuous density, a uniform bound on the density over every $N_{\c
^*}$ provides a local control of $p(t)$ with a polynomial function of
degree 1. This idea is developed in the following example:

\begin{ex}[(Continuous densities, $\mathcal{H}=\mathbb
{R}^d$)]\label{ex3}
Assume that $\mathcal{H} = \mathbb{R}^d$, $P$ has a
continuous density $f$ and is $M$-bounded, and that $\mathcal{M}$ is
finite. In this case, for every $\c^*$, $F_{\c^*}(t) = P(\mathcal
{B}(N_{\c^*},t))$ is differentiable at $0$, with derivative
\[
F_{\c^*}'(0) = \int_{N_{\c^*}}{f(u)\,d
\lambda_{d-1}(u)},
\]
where $ \lambda_{d-1}$ denotes the $(d-1)$-dimensional Lebesgue
measure, considered over the $(d-1)$-dimensional space $N_{\c^*}$.
Therefore, if $P$ satisfies
%
\begin{equation}
\label{continuousdensity} \int_{N_{\c^*}}{f(u)\,d\lambda_{d-1}(u)} <
\frac{B p_{\min}}{128 M^2},
\end{equation}
for every $\c^*$, then there exists $r_0 >0$ such that $P$ satisfies
(\ref{majorationkappa}). It can easily be deduced from (\ref
{continuousdensity}) that a uniform bound on the density located at
$\bigcup_{\c^*} N_{\c^*}$ can provide a sufficient condition for a
distribution $P$ to satisfy a margin condition. Such a result has to be
compared to Theorem~3.2 of \cite{Levrard12}, where it was required
that, for every $\c^*$,
\[
\llVert f_{\mid N_{\c^*} }\rrVert _{\infty} \leq\frac{\Gamma (
d/2  ) B }{2^{d+5} M^{d+1} \pi^{d/2}}p_{\min},
\]
where $\Gamma$ denotes the Gamma function, and $f_{\mid  N_{\c^*}}$ denotes the restriction of $f$ to the set $N_{\c^*}$. Note
however that the uniform bound mentioned above ensures that the Hessian
matrices of the risk function $R$, at optimal codebooks, are positive
definite. This does not necessarily imply that (\ref
{continuousdensity}) is satisfied.

Another interesting parameter of $P$ from the quantization viewpoint is
the following separation factor. It quantifies the difference between
optimal codebooks and local minimizers of the risk.
\end{ex}

\begin{Def}\label{epsilonseparation}
Denote by $\tilde{\mathcal{M}}$ the set of local minimizers of the
map distortion $\c\longmapsto P\gamma(\c,\cdot)$. Let $\varepsilon>0$, then $P$ is
said to be $\varepsilon$-separated if
%
\begin{equation}
\label{separationfactor} \inf_{\c\in\tilde{\mathcal{M}} \cap\mathcal{M}^c}{\ell\bigl(\c,\c ^*\bigr)} =
\varepsilon.
\end{equation}
\end{Def}

It may be noticed that local minimizers of the risk function satisfy
the centroid condition, or have empty cells. Whenever $\mathcal
{H}=\mathbb{R}^d$, $P$ has a density and $P\llVert  x\rrVert  ^2 < \infty$, it can
be proved that the set of minimizers of $R$ coincides with the set of
codebooks satisfying the centroid condition, also called stationary
points (see, e.g., Lemma~A of \cite{Pollard82}). However, this result
cannot be extended to noncontinuous distributions, as proved in
Example 4.11 of \cite{GL00}.

The main results of this paper are based on the following proposition,
which connects the margin condition stated in Definition~\ref
{margincondition} to the previous conditions in \cite{Pollard82} or
\cite{Antos04}. Recall that $k\geq2$.

%
\begin{prop}\label{lienconditionmarge}
Assume that $P$ satisfies a margin condition with radius $r_0$, then
the following properties hold.
\begin{longlist}[(iii)]
\item[(i)] For every $\c^*$ in $\mathcal{M}$ and $\c$ in $\mathcal
{B}(0,M)^k$, if $\llVert   \c- \c^* \rrVert   \leq\frac{Br_0}{4 \sqrt{2}M}$, then
%
\begin{equation}
\label{localconvexity} \ell\bigl(\c,\c^*\bigr) \geq\frac{p_{\min}}{2} \bigl\llVert \c-
\c^*\bigr\rrVert ^2.
\end{equation}

\item[(ii)] $\mathcal{M}$ is finite.
\item[(iii)] There exists $\varepsilon>0$ such that $P$ is
$\varepsilon$-separated.
\item[(iv)] For all $\c$ in $\mathcal{B}(0,M)^k$,
%
\begin{equation}
\label{lienantos} \frac{1}{16M^2} \Var\bigl(\gamma(\c,\cdot) - \gamma\bigl(\c^*(\c),\cdot
\bigr)\bigr) \leq\bigl\llVert \c- \c^*(\c) \bigr\rrVert ^2 \leq
\kappa_0 \ell\bigl(\c,\c^*\bigr),
\end{equation}
where $\kappa_0 = 4kM^2 ( \frac{1}{\varepsilon} \vee\frac{64
M^2}{p_{\min} B^2 r_0^2}  )$, and $\c^*(\c) \in{\arg\min}_{ \c
^* \in\mathcal{M}}{\llVert   \c- \c^*\rrVert  }$.
\end{longlist}
\end{prop}

As a consequence, (\ref{lienantos}) ensures that (\ref
{Antoscondition}) is satisfied, with known constant, which is the
condition required in Theorem~2 of \cite{Antos04}. Moreover, if
$\mathcal{H} = \mathbb{R}^d$, $P$ has a unique optimal codebook up to
relabeling, and has a continuous density, (\ref{localconvexity})
ensures that the second derivative matrix of $R$ at the optimal
codebook is positive definite, with minimum eigenvalue larger than
$p_{\min}/2$. This is the condition required in \cite{Chou94} for
$n\ell(\hat{\c}_n,\c^*)$ to converge in distribution.

It is worth pointing out that the dependency of $\kappa_0$ on
different parameters of $P$ is known. This fact allows us to roughly
discuss how $\kappa_0$ should scale with the parameters $k$, $d$ and
$M$, in the finite-dimensional case. According to Theorem~6.2 of \cite
{GL00}, $R(\c^*)$ scales like $M^2k^{-2/d}$, when $P$ has a density.
Furthermore, it is likely that $r_0 \sim B$ (see,\vspace*{1pt} e.g., the
distributions exposed in Section~\ref{minimaxsection}). Considering
that $\varepsilon\sim R(\c^*) \sim M^2 k^{-2/d}$, $r_0 \sim B \sim M
k^{-1/d}$, and $p_{\min} \sim1/k$ leads to
\[
\kappa_0 \sim k^{2 +4/d}.
\]
At first sight, $\kappa_0$ does not scale with $M$, and seems to
decrease with the dimension, at least in the finite-dimensional case.
However, there is no result on how $\kappa_0$ should scale in the
infinite-dimensional case. Proposition~\ref{lienconditionmarge} allows
us to derive explicit upper bounds on the excess risk in the following section.

\section{Results}\label{Results}

\subsection{Risk bound}\label{riskbound}

The main result of this paper is the following.

\begin{teo}\label{mainresult}
Assume that $k \geq2$, and that $P$ satisfies a margin condition with
radius $r_0$. Let $\kappa_0$ be defined as
\[
\kappa_0 = 4kM^2 \biggl( \frac{1}{\varepsilon} \vee
\frac{64
M^2}{p_{\min} B^2 r_0^2} \biggr).
\]
If $\hat{\c}_n$ is an empirical risk minimizer, then, with
probability larger than $1-e^{-x}$,
%
\begin{equation}
\label{secondinequality} \ell\bigl(\hat{\c}_n,\c^*\bigr) \leq C_0
\kappa_0 \frac{ ( k + \log
 ( \llvert   \bar{\mathcal{M}} \rrvert    )  )M^2}{n} + ( 9 \kappa_0 + 4 )
\frac{16 M^2}{n}x,
\end{equation}
where $C_0$ is an absolute constant.
\end{teo}

This result is in line with Theorem~3.1 in \cite{Levrard12} or Theorem~1 in \cite{Chichi13}, concerning the dependency on the sample size $n$
of the loss $\ell(\hat{\c}_n,\c^*)$. The main advance lies in the
detailed dependency on other parameters of the loss of $\hat{\c}_n$.
This provides a nonasymptotic bound for the excess risk.

To be more precise, Theorem~3.1 in \cite{Levrard12} states that
\[
\mathbb{E}\ell\bigl(\hat{\c}_n,\c^*\bigr) \leq C(k,d,P)M^2/n,
\]
in the finite-dimensional case, for some unknown constant $C(k,d,P)$.
In fact, this result relies on the application of Dudley's entropy
bound. This technique was already the main argument in \cite
{Pollard82} or \cite{Chichi13}, and makes use of covering numbers of
the $d$-dimensional Euclidean unit ball. Consequently, $C(k,d,P)$
strongly depends on the dimension of the underlying Euclidean space in
these previous results. As suggested in \cite{Biau08} or \cite
{Canas12}, the use of metric entropy techniques to derive bounds on the
convergence rate of the distortion may be suboptimal, as it does not
take advantage of the Hilbert space structure of the squared distance
based quantization. This issue can be addressed by using a technique
based on comparison with Gaussian vectors, as done in \cite{Canas12}.
Theorem~\ref{mainresult} is derived that way, providing a
dimension-free upper bound which is valid over separable Hilbert spaces.

It may be noticed that most of results providing slow convergence
rates, such as Theorem~2.1 in \cite{Biau08} or Corollary 1 in \cite
{Linder94}, give bounds on the distortion which do not depend on the
number of optimal codebooks. Theorem~\ref{mainresult} confirms that
$\llvert   \bar{\mathcal{M}} \rrvert  $ is also likely to play a minor
role on the convergence rate of the distortion in the fast rate case.

Another interesting point is that Theorem~\ref{mainresult} does not
require that $P$ has a density or is distributed over points, contrary
to the requirements of the previous bounds in \cite{Pollard82,Antos04} or \cite{Chichi13} which achieved the optimal rate of
$\mathcal{O}(1/n)$. Up to our knowledge, the more general result is to
be found in Theorem~2 of \cite{Antos04}, which derives a convergence
rate of $\mathcal{O}(\log(n)/n)$ without the requirement that $P$ has
a density. It may also be noted that Theorem~\ref{mainresult} does not
require that $\bar{\mathcal{M}}$ contains a single element, contrary
to the results stated in \cite{Pollard82}. According to Proposition
2.2, only (\ref{majorationkappa}) has to be proved for $P$ to satisfy
the assumptions of Theorem~\ref{mainresult}. Since proving that $\llvert  \bar{\mathcal{M}} \rrvert   = 1$ may be difficult, even for simple
distributions, it seems easier to check the assumptions of Theorem~\ref
{mainresult} than the assumptions required in \cite{Pollard82}. An
illustration of this point is given in Section~\ref{Gaussianmixture}.

As will be shown in Proposition~\ref{minimax}, the dependency on
$\varepsilon$ turns out to be sharp when $\varepsilon\sim n^{-1/2}$.
In fact, tuning this separation factor is the core of the demonstration
of the minimax results in \cite{Bartlett98} or \cite{Antos05}.

\subsection{Minimax lower bound}\label{minimaxsection}

This subsection is devoted to obtaining a minimax lower bound on the
excess risk over a set of distributions with continuous densities,
unique optimal codebook, and satisfying a margin condition, in which
some parameters, such as $p_{\min}$ are fixed or uniformly
lower-bounded. It has been already proved in Theorem~4 of \cite
{Antos05} that the minimax distortion over distributions with uniformly
bounded continuous densities, unique optimal codebooks (up to
relabeling), and such that the minimum eigenvalues of the second
derivative matrices at the optimal codebooks are uniformly
lower-bounded, is $\Omega(1/\sqrt{n})$, in the case where $k=3$ and
$d=1$. Extending the distributions used in Theorem~4 of \cite
{Antos05}, Proposition~\ref{minimax} below generalizes this result in
arbitrary dimension $d$, and provides a lower bound over a set of
distributions satisfying a uniform margin condition.

Throughout this subsection, only the case $\mathcal{H} = \mathbb
{R}^d$ is considered, and $\hat{\c}_n$ will denote an empirically
designed codebook, that is a map from $(\mathbb{R}^d)^n$ to $(\mathbb
{R}^d)^k$. Let $k$ be an integer such that $k \geq3$, and $M>0$. For
simplicity, $k$ is assumed to be divisible by $3$. Let us introduce the
following quantities:
\begin{eqnarray*}
&&\cases{\displaystyle m = \frac{2k}{3},
\vspace*{3pt}\cr
\displaystyle\Delta=\frac{5M}{32m^{1/d}}.}
\end{eqnarray*}

To focus on the dependency on the separation factor $\varepsilon$, the
quantities involved in Definition~\ref{margincondition} are fixed as
%
\begin{equation}
\label{conditionmargeminimax}
\cases{\displaystyle B=\Delta,
\vspace*{3pt}\cr
\displaystyle r_0= \frac{7 \Delta}{16},
\vspace*{3pt}\cr
\displaystyle p_{\min} \geq\frac{3}{4k}.}
\end{equation}
Denote by $\mathcal{D}(\varepsilon)$ the set of probability
distributions which are $\varepsilon$-separated, have continuous
densities and unique optimal codebooks, and which satisfy a margin
condition with parameters defined in (\ref{conditionmargeminimax}).
The minimax result is the following.

\begin{prop}\label{minimax}
Assume that $k \geq3$ and $n \geq3k/2$. Then, for any empirically
designed codebook,
\[
\sup_{P \in\mathcal{D}(c_1/\sqrt{n})} \mathbb{E} \ell\bigl(\hat{\c }_n,\c^*
\bigr) \geq c_0 M^2 \frac{\sqrt{k^{1-{4/d}}}}{\sqrt{n}},
\]
where $c_0 >0 $ is an absolute constant, and
\[
c_1 = \frac{(5M)^2}{4(32m^{1/4+1/d})^2}.
\]
\end{prop}

Proposition~\ref{minimax} is in line with the previous minimax lower
bounds obtained in Theorem~1 of \cite{Bartlett98} or Theorem~4 of
\cite{Antos05}. Proposition~\ref{minimax}, as well as these two
previous results, emphasizes the fact that fixing the parameters of the
margin condition uniformly over a class of distributions does not
guarantee an optimal uniform convergence rate. This shows that a
uniform separation assumption is needed to derive a sharp uniform
convergence rate over a set of distributions.

Furthermore, as mentioned above, Proposition~\ref{minimax} also
confirms that the minimax distortion rate over the set of distributions
with continuous densities, unique optimal codebooks, and such that the
minimum eigenvalues of the Hessian matrices are uniformly lower bounded
by $3/8k$, is still $\Omega(1/\sqrt{n})$ in the case where $d>1$ and
$k\geq3$.

This minimax lower bound has to be compared to the upper risk bound
obtained in Theorem~\ref{mainresult} for the empirical risk minimizer
$\hat{\c}_n$, over the set of distributions $\mathcal{D}(c_1/\sqrt
{n})$. To be more precise, Theorem~\ref{mainresult} ensures that,
provided that $n$ is large enough,
\[
\sup_{P \in\mathcal{D}(c_1/\sqrt{n})}\mathbb{E} \ell\bigl(\hat{\c }_n,\c^*
\bigr) \leq\frac{g(k,d,M)}{\sqrt{n}},
\]
where $g(k,d,M)$ depends only on $k$, $d$ and $M$. In other words, the
dependency of the upper\vspace*{1pt} bounds stated in Theorem~\ref{mainresult} on
$\varepsilon$ turns out to be sharp whenever $\varepsilon\sim
n^{-1/2}$. Unfortunately, Proposition~\ref{minimax} cannot be
easily extended to the case where $\varepsilon\sim n^{-\alpha}$, with
$0<\alpha<1/2$. Consequently, an open question is whether the upper
bounds stated in Theorem~\ref{mainresult} remains accurate with
respect to $\varepsilon$ in this case.

\subsection{Quasi-Gaussian mixture example}\label{Gaussianmixture}

The aim of this subsection is to illustrate the results exposed in
Section~\ref{Results} with Gaussian mixtures in dimension $d=2$. The
Gaussian mixture model is a typical and well-defined clustering example.

In general, a Gaussian mixture distribution $\tilde{P}$ is defined by
its density
\[
\tilde{f}(x) = \sum_{i=1}^{\tilde{k}}{
\frac{\theta_i}{2 \pi\sqrt
{ \llvert   \Sigma_i \rrvert   }}e^{-(1/2)(x-m_i)^t \Sigma_i^{-1}
(x-m_i)}},
\]
where $\tilde{k}$ denotes the number of components of the mixture, and
the $\theta_i$'s denote the weights of the mixture, which satisfy
$\sum_{i=1}^{k}{\theta_i} = 1$. Moreover, the $m_i$'s denote the
means of the mixture, so that $m_i \in  \mathbb{R}^2$, and the
$\Sigma_i$'s are the $2\times2$ variance matrices of the components.

We restrict ourselves to the case where the number of components
$\tilde{k}$ is known, and match the size $k$ of the codebooks. To ease
the calculation, we make the additional assumption that every component
has the same diagonal variance matrix $\Sigma_i = \sigma^2 I_2$. Note
that a similar result to Proposition~\ref{Gaussianboundary} can be
derived for distributions with different variance matrices $\Sigma_i$,
at the cost of more computing.

Since the support of a Gaussian random variable is not bounded, we
define the ``quasi-Gaussian'' mixture model as follows, truncating each
Gaussian component.
Let the density $f$ of the distribution $P$ be defined by
\[
f(x) = \sum_{i=1}^{k}{\frac{\theta_i}{ 2 \pi\sigma^2 N_i}
e^{-{\llVert  x-m_i\rrVert  ^2}/(2 \sigma^2)}}\mathbh{1}_{\mathcal{B}(0,M)},
\]
where $N_i$ denotes a normalization constant for each Gaussian variable.

Let $\eta$ be defined as $\eta= 1 - \min_{i=1, \ldots, k}{N_i}$.
Roughly, the model proposed above will be close the Gaussian
mixture model when $\eta$ is small.
Denote by $\tilde{B} = {\inf_{i \neq j}}{\llVert  m_i - m_j\rrVert  }$ the smallest
possible distance between two different means of the mixture. To avoid
boundary issues we assume that, for all $i = 1, \ldots, k$, $\mathcal
{B}(m_i, \tilde{B}/3) \subset\mathcal{B}(0,M)$.

Note that the assumption $\mathcal{B}(m_i, \tilde{B}/3) \subset
\mathcal{B}(0,M)$ can easily be satisfied if $M$ is chosen large
enough. For such a model, Proposition~\ref{Gaussianboundary} offers a
sufficient condition for $P$ to satisfy a margin condition.

%
\begin{prop}\label{Gaussianboundary}
Let $\theta_{\min} = \min_{i=1,\ldots,k}{\theta_i}$, and $\theta
_{\max}=\max_{i=1,\ldots,k}{\theta_i}$. Assume that
%
\begin{equation}
\label{gaussiancondition} \frac{\theta_{\min}}{\theta_{\max}} \geq\frac{2048 k}{(1-\eta)
\tilde{B}} \max{ \biggl(
\frac{ \sigma^2}{ \tilde{B}(1 - e^{-\tilde
{B}^2/{2048\sigma^2}})}, \frac{ k M^3}{7\sigma^2 (e^{\tilde
{B}^2/{32\sigma^2}}-1)} \biggr) }.
\end{equation}
Then $P$ satisfies a margin condition with radius $\frac{\tilde{B}}{8}$.
\end{prop}

It is worth mentioning that $P$ has a continuous density, and that
according to~(i) in Proposition~\ref{lienconditionmarge}, the second
derivative matrices of the risk function, at the optimal codebooks,
must be positive definite. Thus, $P$ might be in the scope of the
result in \cite{Pollard82}. However, there is no elementary proof of
the fact that $ \llvert   \bar{\mathcal{M}} \rrvert  =1$, whereas $\mathcal{M}$ is
finite is guaranteed by Proposition~\ref{lienconditionmarge}. This
shows that the margin condition given in Definition~\ref
{margincondition} may be easier to check than the condition presented
in \cite{Pollard82}. The condition (\ref{gaussiancondition}) can be
decomposed as follows. If
\[
\frac{\theta_{\min}}{\theta_{\max}} \geq\frac{2048 k \sigma
^2}{(1-\eta) \tilde{B}^2(1 - e^{-\tilde{B}^2/{2048\sigma^2}})},
\]
then every optimal codebook $\c^*$ must be close to the vector of
means of the mixture $\mathbf{m} = (m_1, \ldots, m_k)$. Therefore, it
is possible to approximately locate the $N_{\c^*}$'s, and to derive an
upper bound on the weight function $p(t)$ defined above Definition~\ref
{margincondition}. This leads to the second term of the maximum in
(\ref{gaussiancondition}).

This condition can be interpreted as a condition on the polarization of
the mixture. A favorable case for vector quantization seems to be when
the poles of the mixtures are well separated, which is equivalent to
$\sigma$ is small compared to $\tilde{B}$, when considering Gaussian
mixtures. Proposition~\ref{Gaussianboundary} gives details on how
$\sigma$ has to be small compared to $\tilde{B}$, in order to satisfy
the requirements of Definition~\ref{margincondition}.

It may be noticed that Proposition~\ref{Gaussianboundary} offers
almost the same condition as Proposition 4.2 in \cite{Levrard12}. In
fact, since the Gaussian mixture distributions have a continuous
density, making use of (\ref{continuousdensity}) in Example~\ref{ex3} ensures
that the margin condition for Gaussian mixtures is equivalent to a
bound on the density over $ \bigcup_{\c^*}N_{\c^*}$.

It is important to note that this result is valid when $k$ is known and
matches exactly the number of components of the mixture. When the
number of code points $k$ is different from the number of components
$\tilde{k}$ of the mixture, we have no general idea of where the
optimal code points can be located.

Moreover, suppose that there exists only one optimal codebook $\c^*$,
up to relabeling, and that we are able to locate this optimal codebook
$\c^*$. As stated in Proposition~\ref{lienconditionmarge}, the key
quantity is in fact $B  = \inf_{i\neq j}\llVert  c^*_i - c^*_j\rrVert  $. In the
case where $\tilde{k} \neq k$, there is no simple relation between
$\tilde{B}$ and $B$. Consequently, a condition like in Proposition
\ref{Gaussianboundary} could not involve the natural parameter of the
mixture $\tilde{B}$.

\section{Proofs}\label{Proofs}
\subsection{Proof of Proposition \protect\ref{minorationBetpmin}}
The lower bound on $B$ follows from a compactness argument for the weak
topology on $\mathcal{H}$, exposed in the following lemma. For the
sake of completeness, it is recalled that a sequence $c_n$ of elements
in $\mathcal{H}$ weakly converges to $c$, denoted by $c_n
\rightharpoonup_{n \rightarrow\infty} c$, if, for every continuous
linear real-valued function $f$, $f(c_n) \rightarrow_{n\rightarrow
\infty} f(c)$. Moreover, a function $\phi$ from $\mathcal{H}$ to
$\mathbb{R}$ is weakly lower semi-continuous if, for all $\lambda \in \mathbb{R}$, the level sets $\{ c \in\mathcal{H}\mid  \phi
(c) \leq\lambda\}$ are closed for the weak topology.

\begin{lem}\label{convergencefaible}
Let $\mathcal{H}$ be a separable Hilbert space, and assume that $P$ is
\mbox{$M$-}bounded. Then:
\begin{longlist}[(iii)]
\item[(i)] $\mathcal{B}(0,R)^k$ is weakly compact, for every $R \geq0$,
\item[(ii)] $\c\mapsto P\gamma(\c,\cdot)$ is weakly lower semi-continuous,
\item[(iii)] $\mathcal{M}$ is weakly compact.
\end{longlist}
\end{lem}

A more general statement of Lemma~\ref{convergencefaible} can be found
in Section~5.2 of \cite{Fischer10}, for quantization with Bregman
divergences. However, since the proof is much simpler in the special
case of the squared-norm based quantization in a Hilbert space, it is
briefly recalled in Section~A.1 (supplementary material \cite{supple}).

Let $\c'_n$ be a sequence of optimal codebooks such that $\llVert   c'_{1,n}
- c'_{2,n}\rrVert   \rightarrow B$, as $n \rightarrow \infty$.
Then, according to Lemma~\ref{convergencefaible}, there exists a
subsequence $\c_n$ and an optimal codebook $\c^*$, such that $\c_n {\rightharpoonup}_{ n \rightarrow\infty} \c^*$, for the weak
topology. Then it is clear that $(c_{1,n} - c_{2,n}) {\rightharpoonup
}_{ n \rightarrow\infty} (c_1^* - c_2^*)$.

Since $u \mapsto\llVert  u\rrVert  $ is weakly lower semi-continuous on $\mathcal
{H}$ (see, e.g., Proposition~3.13 in \cite{Brezis11}), it follows that
\[
\bigl\llVert c_1^* - c_2^* \bigr\rrVert \leq\underset{ n
\rightarrow\infty} {\lim\inf} \llVert c_{1,n} - c_{2,n} \rrVert
= B.
\]
Noting that $\c^*$ is an optimal codebook, and the support of $P$ has
more than $k$ points, Proposition 1 of \cite{Graf07} ensures that $\llVert
c_1^* - c_2^*\rrVert   >0 $.\vspace*{1pt}

The uniform lower bound on $p_{\min}$ follows from the argument that,
since the support of $P$ contains more than $k$ points, then $R^*_k <
R^*_{k-1}$, where\vspace*{2pt} $R^*_j$ denotes the minimum distortion achievable for
$j$-points quantizers (see, e.g., Proposition~1 in \cite{Graf07}).
Denote by $\alpha$ the quantity $ R^*_{k-1}-R^*_k$, and suppose that
$p_{\min} < \frac{\alpha}{4M^2}$. Then there exists an optimal
codebook of size $k$, $\c^{*,k}= (c^{*,k}_1, \ldots, c^{*,k}_k)$,
such that $ P(V_1(\c^{*,k})) < \frac{\alpha}{4M^2} $. Let\vspace*{2pt} $\c
^{*,{k-1}}$ denote an optimal codebook of size $(k-1)$, and define the
following $k$-points quantizer:
\begin{eqnarray*}
\cases{\displaystyle Q(x)= c_1^{*,k}, &\quad if $x \in V_1
\bigl(\c^{*,k}\bigr)$,
\vspace*{3pt}\cr
\displaystyle Q(x)= c_j^{*,k-1}, &\quad
if $x \in V_j\bigl(\c^{*,{k-1}}\bigr) \cap \bigl(
V_1\bigl(\c^{*,k}\bigr) \bigr)^{c}$.}
\end{eqnarray*}
Since $P(\partial V_1(\c^{*,k}))= P(\partial V_j(\c^{*,{k-1}})) = 0$,
for $j=1, \ldots, k-1$, Q is defined $P$ almost surely. Then it is
easy to see that
\[
R(Q) \leq P\bigl(V_1\bigl(\c^{*,k}\bigr)\bigr) 4
M^2 + R^*_{k-1} < R^*_{k}.
\]
Hence, the contradiction. Therefore, we have $p_{\min} \geq\frac
{\alpha}{4M^2}$.

\subsection{Proof of Proposition \protect\ref{lienconditionmarge}}\label{ProofofProposition{lienconditionmarge}}

The proof of (i) in Proposition~\ref{lienconditionmarge} is based on
the following lemma.

%
\begin{lem}\label{boundarycloseness}
Let $\c$ and $\c^*$ be in $\mathcal{B}(0,M)^k$, and $x \in V_i(\c
^*) \cap V_j(\c) \cap\mathcal{B}(0,M)$, for $i \neq j$. Then
%
\begin{eqnarray}
\biggl\llvert \biggl\langle x - \frac{c_i + c_j}{2}, c_i -
c_j \biggr\rangle \biggr\rrvert & \leq&4 \sqrt{2}M \bigl\llVert \c-
\c^* \bigr\rrVert \label{Vor1},
\\
d\bigl(x, \partial V_i\bigl(\c^*\bigr)\bigr) & \leq&
\frac{4 \sqrt{2} M}{B} \bigl\llVert \c- \c ^* \bigr\rrVert \label{Vor2}.
\end{eqnarray}
\end{lem}

The two statements of Lemma~\ref{boundarycloseness} emphasize the fact
that, provided that $\c$~and~$\c^*$ are quite similar, the areas on
which the labels may differ with respect to $\c$ and $\c^*$ should be
close to the boundary of Voronoi diagrams. This idea is mentioned in
the proof of Corollary 1 in \cite{Antos04}. Nevertheless, we provide a
simpler proof in Section~A.2 (supplementary material \cite{supple}).

Equipped\vspace*{1pt} with Lemma~\ref{boundarycloseness}, we are in a position to
prove (\ref{localconvexity}).
Let $\c$ be in $\mathcal{B}(0,M)^k$, and $(W_1(\c), \ldots, W_k(\c
))$ be a Voronoi partition associated with $\c$, as defined in
Section~\ref{Notation}. Let $\c^*$ be in $\mathcal{M}$, then $\ell
(\c,\c^*)$ can be decomposed as follows:
\begin{eqnarray*}
P \gamma(\c,\cdot) &=&\sum_{i=1}^{k}{P\bigl(
\llVert x - c_i \rrVert ^2 \mathbh {1}_{W_i(\c)}(x)
\bigr)}
\\
& =& \sum_{i=1}^{k}{P\bigl(\llVert x -
c_i \rrVert ^2 \mathbh{1}_{V_i(\c^*)}(x)\bigr)}
\\
&&{}+ \sum_{i=1}^{k}{P\bigl(\llVert
x-c_i\rrVert ^2\bigl(\mathbh{1}_{W_i(\c)}(x)-
\mathbh {1}_{V_i(\c^*)}(x)\bigr)\bigr)}.
\end{eqnarray*}
Since, for all $i = 1,\ldots, k$, $P(x\mathbh{1}_{V_i(\c^*)}(x)) =
P(V_i(\c^*))c_i^*$ (centroid condition), we may write
\begin{eqnarray*}
&& P\bigl(\llVert x - c_i \rrVert ^2 \mathbh{1}_{V_i(\c^*)}(x)
\bigr)
\\
&&\qquad = P\bigl(V_i\bigl(\c^*\bigr)\bigr) \bigl\llVert
c_i - c_i^*\bigr\rrVert ^2 + P\bigl( \bigl
\llVert x - c_i^* \bigr\rrVert ^2 \mathbh{1}_{V_i(\c^*)}(x)
\bigr),
\end{eqnarray*}
from which we deduce that
\begin{eqnarray*}
P \gamma(\c,\cdot) &=& P \gamma\bigl(\c^*,\cdot\bigr) + \sum_{i=1}^{k}{P
\bigl(V_i\bigl(\c^*\bigr)\bigr) \bigl\llVert c_i -
c_i^*\bigr\rrVert ^2}
\\
&&{}+ \sum_{i=1}^{k}{P\bigl(\llVert
x-c_i\rrVert ^2\bigl(\mathbh{1}_{W_i(\c
)}(x)-
\mathbh{1}_{V_i(\c^*)}(x)\bigr)\bigr)},
\end{eqnarray*}
which leads to
\begin{eqnarray*}
\ell\bigl(\c,\c^*\bigr) &\geq& p_{\min} \bigl\llVert \c- \c^*\bigr\rrVert
^2
\\
&&{}+ \sum_{i=1}^{k}{ \sum
_{j \neq i} {P \bigl(\bigl(\llVert x-c_j\rrVert
^2 - \llVert x-c_i\rrVert ^2\bigr)\mathbh
{1}_{V_i(\c^*)\cap W_j(\c)}(x) \bigr)}}.
\end{eqnarray*}
Since $x\in W_j(\c) \subset V_j(\c)$, $\llVert  x-c_j\rrVert  ^2 - \llVert  x-c_i\rrVert  ^2 \leq
0$. Thus, it remains to bound from above
\[
\sum_{i=1}^{k}{ \sum
_{j \neq i} {P \bigl(\bigl(\llVert x-c_i\rrVert
^2 - \llVert x-c_j\rrVert ^2\bigr)
\mathbh{1}_{V_i(\c^*)\cap W_j(\c)}(x) \bigr)}}.
\]
Noticing that
\[
\llVert x-c_i\rrVert ^2 - \llVert x-c_j
\rrVert ^2 = 2 \biggl\langle c_j -c_i, x -
\frac
{c_i+c_j}{2} \biggr\rangle,
\]
and using Lemma~\ref{boundarycloseness}, we get
\begin{eqnarray*}
&& \bigl(\llVert x-c_i\rrVert ^2 - \llVert
x-c_j\rrVert ^2\bigr)\mathbh{1}_{V_i(\c^*)\cap W_j(\c)}(x)
\\
&&\qquad \leq8 \sqrt{2}M \bigl\llVert \c- \c^* \bigr\rrVert \mathbh{1}_{ V_i(\c^*)\cap W_j(\c
) \cap N_{\c^*} ( (4\sqrt{2}M/B) \llVert   \c- \c^* \rrVert
) }(x).
\end{eqnarray*}
Hence,
\begin{eqnarray*}
&& \sum_{i=1}^{k}{P \bigl( \llVert
x-c_i\rrVert ^2\bigl(\mathbh{1}_{W_i(\c)}(x)-
\mathbh {1}_{V_i(\c^*)}(x)\bigr) \bigr)}
\\
&&\qquad \geq- 8 \sqrt{2}M \bigl\llVert \c- \c^* \bigr\rrVert p \biggl(
\frac{4 \sqrt{2}M}{B}\bigl\llVert \c- \c^* \bigr\rrVert \biggr).
\end{eqnarray*}
Consequently, if $P$ satisfies (\ref{majorationkappa}), then, if $\llVert
\c- \c^* \rrVert   \leq\frac{B r_0}{4\sqrt{2}M}$, it follows that
\begin{eqnarray*}
&&\ell\bigl(\c,\c^*\bigr) \geq\frac{p_{\min}}{2}\bigl\llVert \c- \c^*\bigr\rrVert
^2,
\end{eqnarray*}
which proves (i).

Suppose that $\mathcal{M}$ is not finite. According to Lemma~\ref
{convergencefaible}, there exists a sequence $\c_n$ of optimal
codebooks and an optimal codebook $\c^*$ such that for all~$n$, $\c_n
\neq\c^*$ and $\c_n \rightharpoonup_{n \rightarrow\infty} \c^*$.
Assume that there exists $i$ in $ \{1, \ldots, k  \}$ such
that $\lim\inf_n \llVert   c_{n,i} \rrVert  ^2 > \llVert   c_i \rrVert  ^2$. Then\vspace*{1pt} $\lim\inf_n
\llVert   x - c_{n,i} \rrVert  ^2 > \llVert   x-c_i \rrVert  ^2$, for every $x$ in $\mathcal{H}$.
Let $x$ be in $\displaystyle\Vo_i(\c)$, and $j \neq i$, then
\[
\underset{ n \rightarrow\infty} {\lim\inf}\llVert x - c_{n,j} \rrVert
^2 \geq \llVert x - c_j \rrVert ^2 > \llVert
x - c_i \rrVert ^2,
\]
which leads to $\lim\inf_n \gamma(\c_n,x) > \gamma(\c,x)$. Since
$P(\displaystyle{\Vo}_i(\c)) > 0$, it easily follows that
\[
\underset{ n \rightarrow\infty} {\lim\inf} P \gamma(\c_n,\cdot) \geq P
\underset{ n \rightarrow\infty} {\lim\inf} \gamma(\c_n,\cdot) > P \gamma(
\c,\cdot),
\]
which is impossible. Hence, there exists a subsequence $\bar{\c}_n$
of $\c_n$ such that, for $i=1, \ldots, k$, $\llVert   \bar{c}_{n,i} \rrVert
\rightarrow_{n \rightarrow\infty} \llVert   c^*_i \rrVert  $. Since Hilbert spaces
are uniformly convex spaces, hence satisfy the Radon--Riesz property
(see, e.g., Propositions 5.1~and~3.32 in~\cite{Brezis11}),
it follows that $\bar{\c}_n \rightarrow_{n \rightarrow\infty} \c
^*$. This contradicts (\ref{localconvexity}) and proves (ii).

The proof of (iii) is based on the following two lemmas.

%
\begin{lem}\label{nocentroid}
Let $\c$ be in ${\mathcal{B} ( \mathcal{M},\frac{B r_0}{4
\sqrt{2}M}  )}$. If $\c$ satisfies the centroid condition, then
$\c$ is in $\mathcal{M}$.
\end{lem}

Lemma~\ref{nocentroid} ensures that no local minimizer with nonempty
cells can be found in a neighborhood of $\mathcal{M}$. We postpone its
proof to Section~A.3 (supplementary material \cite{supple}). Lemma~\ref{localminimum} below shows that
the infimum distortion over codebooks which are away from $\mathcal
{M}$ is achieved.

%
\begin{lem}\label{localminimum}
For every $r>0$, there exists $\c_r$ in $\mathcal{B}(0,M+r)^k
\setminus\mathcal{B}^o(\mathcal{M},r)$ such that
\[
\inf_{\mathcal{H}^k \setminus\mathcal{B}^o(\mathcal{M},r)}{P\gamma (\c,\cdot)} = P \gamma(\c_r,\cdot).
\]
\end{lem}

The proof of Lemma~\ref{localminimum} is given in Section~A.4 (supplementary material \cite{supple}).
Let $\tilde{\c} \notin \mathcal{M}$ be a local minimizer of the
distortion. If $\tilde{\c}$ has empty cells, then $P \gamma(\tilde
{\c},\cdot) \geq R^*_{k-1} > R^*_{k}$. Assume that $\tilde{\c}$ has no
empty cells. Then $\tilde{\c}$ satisfies the centroid condition, thus
Lemma~\ref{nocentroid} ensures that $\llVert   \tilde{\c} - \c^* \rrVert   \geq
r$, for every optimal codebook $\c^*$ and for $r=\frac{B r_0}{4 \sqrt
{2}M}$. Lemma~\ref{localminimum} provides $\c_r$ such that $P \gamma
(\tilde{\c},\cdot) \geq P\gamma(\c_r,\cdot)>0$. Hence, (iii) is proved.

The left part of (\ref{lienantos}) follows from the elementary inequality
%
\begin{equation}
\label{majorationdebase} \forall x \in\mathcal{B}(0,M) \qquad\bigl\llvert \gamma(\c,x) -
\gamma \bigl(\c^*(\c),x\bigr) \bigr\rrvert \leq4M \max_{i=1, \ldots, k}{
\bigl\llVert c_i - c_i^*(\c) \bigr\rrVert }.
\end{equation}

According\vspace*{2pt} to (\ref{localconvexity}), if $\llVert  \c- \c^*(\c)\rrVert   \leq
\frac{B r_0}{4 \sqrt{2}M}$, then $\ell(\c,\c^*) \geq\frac
{p_{\min}}{2}\llVert   \c- \c^*(\c) \rrVert  ^2$.
Now turn to the case where $\llVert   \c- \c^*(\c) \rrVert   \geq\frac{B
r_0}{4\sqrt{2}M}=r$. Then\vspace*{2pt} Lemma~\ref{localminimum} provides $\c_r$
such that $\ell(\c,\c^*) \geq\ell(\c_r,\c^*)$. Such a $\c_r$ is
a local minimum of $\c\mapsto P\gamma(\c,\cdot)$, or satisfies $\llVert  \c_r
- \c^*(\c_r)\rrVert  =r$. Hence, we deduce
\begin{eqnarray*}
\ell\bigl(\c,\c^*\bigr) \geq\ell\bigl(
\c_r,\c^*\bigr)& \geq&\varepsilon\wedge\frac
{p_{\min}}{2}r^2
 \geq \biggl(\varepsilon\wedge\frac{p_{\min} B^2 r_0 ^2}{64 M^2} \biggr) \frac{\llVert   \c- \c^*(\c) \rrVert  ^2}{4kM^2}.
\end{eqnarray*}
Note that, since $B \leq2M$ and $r_0 \leq2M$, $ (\varepsilon
\wedge\frac{p_{\min} B^2 r_0 ^2}{64 M^2}  )/4kM^2 \leq
p_{\min}/2$. This proves (\ref{lienantos}).

\subsection{Proof of Theorem~\protect\ref{mainresult}}

Throughout this subsection, $P$ is assumed to satisfy a margin
condition with radius $r_0$, and we\vspace*{1pt} denote by $\varepsilon$ its
separation factor.
A~nondecreasing map $\Phi\dvtx \mathbb{R}^+ \rightarrow\mathbb{R}^+$ is
called sub-root if
$x \mapsto\frac{\Phi(x)}{\sqrt{x}}$ is nonincreasing.

The following localization theorem, derived from Theorem~6.1 in \cite
{Blanchard08}, is the main argument of our proof.

\begin{teo}\label{localization}
Let $\mathcal{F}$ be a class of uniformly bounded measurable functions
such that there exists $\omega\dvtx  \mathcal{F} \longrightarrow\mathbb
{R}^+$ satisfying
\[
\forall f \in\mathcal{F}, \qquad\Var(f) \leq\omega(f).
\]
Assume that
\[
\forall r >0, \qquad\mathbb{E} \Bigl( \sup_{\omega(f)\leq r}{\bigl\llvert
(P-P_n)f \bigr\rrvert } \Bigr) \leq\Phi(r),
\]
for some sub-root function $\Phi$. Let $K$ be a positive constant, and
denote by $r^*$ the unique solution of the equation $\Phi(r) = r/24K$.

Then, for all $x>0$, with probability larger than $1-e^{-x}$,
\[
\forall f \in\mathcal{F}, \qquad Pf - P_n f \leq K^{-1}
\biggl( \omega (f) + r^* + \frac{(9K^2+16K \sup_{f \in\mathcal{F}}{\llVert  f\rrVert  _{\infty
})x}}{4n} \biggr).
\]
\end{teo}

A proof of Theorem~\ref{localization} is given in Section~5.3 of \cite
{Levrard12}. The proof of (\ref{secondinequality}) follows from the
combination of Proposition~\ref{lienconditionmarge} and a direct
application of Theorem~\ref{localization}. To be more precise, let
$\mathcal{F}$ denote the set
\[
\mathcal{F} = \bigl\{ \gamma(\c,\cdot) - \gamma\bigl(\c^*(\c),\cdot\bigr) \mid \c
\in\mathcal{B}(0,M)^k \bigr\}.
\]
According to (\ref{majorationdebase}), it is clear that, for every $f \in \mathcal{F}$,
\[
\cases{\displaystyle \llVert f\rrVert _{\infty} \leq 8 M^2,
\vspace*{3pt}\cr
\Var(f)
\leq 16 M^2 \bigl\llVert \c- \c^*(\c) \bigr\rrVert ^2.}
\]
Define $\omega(f) = 16 M^2 \llVert   \c- \c^*(\c) \rrVert  ^2$. It remains to
bound from above the complexity term. This is done in the following proposition.

\begin{prop}\label{chainagegaussien}
One has
%
\begin{equation}
\label{definitiondelta} \mathbb{E} \sup_{f \in\mathcal{F}, \omega(f) \leq\delta}{\bigl\llvert
(P-P_n) f \bigr\rrvert } \leq\frac{  (4 \sqrt{\pi k} + \sqrt{2\log ( \llvert  \bar
{\mathcal{M}} \rrvert    ) }  )}{\sqrt{n}}{\sqrt{\delta}}.
\end{equation}
\end{prop}

The proof of Proposition~\ref{chainagegaussien} relies on the use of
Gaussian complexities combined with Slepian's lemma (see, e.g., Theorem~3.14 in \cite{Massart03}), as done in \cite{Canas12}. We postpone it
to the following subsection. Let $\Phi$ be defined as the right-hand
side of (\ref{definitiondelta}), and let $\delta^*$ denote the
solution of the equation $\Phi(\delta) = \delta/ 24K$, for some
positive $K >0$. Then $\delta^*$ can be expressed as
\[
\delta= \frac{576 K^2}{n} \bigl(4 \sqrt{\pi k} + \sqrt{2\log \bigl( \llvert
\bar{\mathcal{M}} \rrvert \bigr)} \bigr)^2 \leq C \frac
{K^2  ( k + \log (\llvert  \bar{\mathcal{M}} \rrvert
)  )}{n}:=
\frac{K^2 \Xi}{n},
\]
where $C=18{,}432 \pi$, and $\Xi=C  ( k + \log (\llvert  \bar
{\mathcal{M}} \rrvert   )  )$. Applying Theorem~\ref
{localization} to $\mathcal{F}$ leads to, with probability larger than
$1-e^{-x}$,
\begin{eqnarray*}
&& (P-P_n) \bigl(\gamma(\mathbf{c},\cdot) - \gamma\bigl(\mathbf{c}^*(\c),\cdot
\bigr)\bigr)
\\
&&\qquad  \leq K^{-1} 16 M^2 \bigl\llVert \mathbf{c}-
\mathbf{c}^*(\c)\bigr\rrVert ^2
+ \frac{ K \Xi
}{n} + \frac{9K+128M^2}{4n}x.
\end{eqnarray*}
Introducing the inequality $ \kappa_0 \ell(\c,\c^*) \geq\llVert   \c- \c
^*(\c) \rrVert  ^2$ provided by Proposition~\ref{lienconditionmarge}, and
choosing $K = 32 M^2 \kappa_0$ leads to (\ref{secondinequality}).

\subsubsection{Proof of Proposition \protect\ref{chainagegaussien}}

As mentioned above, this proof relies on the use of Gaussian
complexities (see, e.g., \cite{Bartlett02}). As will be shown below,
avoiding Dudley's entropy argument by introducing some Gaussian random
vectors allows us to take advantage of the underlying Hilbert space
structure. The first step is to decompose the complexity term according
to optimal codebooks in the following way:
\begin{eqnarray*}
&& \mathbb{E} \sup_{\llVert   \c- \c^*(\c)\rrVert  ^2 \leq\delta/16 M^2}{\bigl\llvert (P-P_n) \bigl(
\gamma(\c,\cdot) - \gamma\bigl(\c^*(\c),\cdot\bigr)\bigr) \bigr\rrvert }
\\
&&\qquad \leq {\mathbb{E}\sup_{\c^* \in\bar{\mathcal{M}}} \sup_{\llVert   \c- \c^*\rrVert
^2 \leq\delta/16 M^2}{
\bigl\llvert (P-P_n) \bigl(\gamma(\c,\cdot) - \gamma\bigl(\c ^*,\cdot\bigr)
\bigr) \bigr\rrvert }}.
\end{eqnarray*}
Let $Y_\c^*$ denote the random variable defined by
\[
Y_\c^* = \sup_{\llVert   \c- \c^*\rrVert  ^2 \leq\delta/16 M^2}{\bigl\llvert
(P-P_n) \bigl(\gamma(\c,\cdot) - \gamma\bigl(\c^*,\cdot\bigr)\bigr) \bigr
\rrvert },
\]
for every $\c^*$ in $\bar{\mathcal{M}}$. It easily follows that
%
\begin{equation}
\label{decompositionwrtM} \mathbb{E} \sup_{\c^* \in\bar{\mathcal{M}}} {Y_\c^*}
\leq \mathbb{E} \sup_{\c^* \in\bar{\mathcal{M}}} { \bigl(Y_\c ^*-
\mathbb{E} Y_\c^* \bigr)} + \sup_{\c^* \in\bar{\mathcal
{M}}}{\mathbb{E}
Y_\c^* }.
\end{equation}
Since, for a fixed $\c^*$, $ \llVert   \gamma(\c,\cdot) - \gamma(\c
^*,\cdot) \rrVert  _{\infty} \leq\sqrt{\delta}$ when $\llVert   \c- \c^* \rrVert
^2 \leq\delta/16M^2$, the bounded difference inequality (see, e.g.,
Theorem~5.1 in \cite{Massart03}) ensures that $Y_\c^*$ is a
sub-Gaussian random variable, with variance bounded from above by
$\delta/n$, that is,
\begin{eqnarray*}
&&\cases{ \displaystyle \mathbb{P} \biggl( Y_\c^* - \mathbb{E} Y_\c^*
\geq\sqrt{\frac
{2\delta x }{n}} \biggr) \leq e^{-x},
\vspace*{3pt}\cr
\displaystyle \mathbb{P} \biggl(
\mathbb{E} Y_\c^* - Y_\c^* \geq\sqrt{
\frac
{2\delta x }{n}} \biggr) \leq e^{-x},}
\end{eqnarray*}
for every $\c^*$ in $\bar{\mathcal{M}}$ and every positive $x$. For
a more general definition of sub-Gaussian random variables, the
interested reader is referred to \cite{Massart03}. Applying Lemma 6.3
in \cite{Massart03} to the special case of sub-Gaussian random
variables leads to
%
\begin{equation}
\label{deviationovercodebooks} \mathbb{E} \sup_{\c^* \in\bar{\mathcal{M}}} { \bigl(Y_\c
^*-\mathbb{E} Y_\c^* \bigr)} \leq\sqrt{\frac{2 \log ( \llvert   \bar{\mathcal{M}} \rrvert    ) \delta}{n}}.
\end{equation}

Next, we bound from above the quantities $\mathbb{E} Y_\c^*$. Let $\c
^*$ be fixed, and let $\sigma_1, \ldots, \sigma_n$ denote some
independent Rademacher variables. According to the symmetrization
principle (see, e.g., Section~2.2 of \cite{Koltchinskii04}),
\begin{eqnarray*}
&& \mathbb{E} \sup_{\llVert   \c- \c^*\rrVert  ^2 \leq\delta/16 M^2}{\bigl\llvert (P-P_n) \bigl(
\gamma(\c,\cdot) - \gamma\bigl(\c^*,\cdot\bigr)\bigr) \bigr\rrvert }
\\
&&\qquad \leq 2 \mathbb{E}_{X,\sigma} \sup_{\llVert   \c- \c^* \rrVert  ^2 \leq\delta/16
M^2}{
\frac{1}{n} \sum_{i=1}^{n}{
\sigma_i \bigl(\gamma(\c,X_i) - \gamma \bigl(
\c^*,X_i\bigr)\bigr)}},
\end{eqnarray*}
where $\mathbb{E}_{Y}$ denotes integration with respect to the
distribution of $Y$. Let $g_1, \ldots, g_n$ denote some independent
standard Gaussian variables. Applying Lemma 4.5 in~\cite{Ledoux91}
leads to
\begin{eqnarray*}
&& \mathbb{E}_{X,\sigma} \sup_{\llVert   \c- \c^* \rrVert  ^2 \leq\delta/16
M^2}{\frac{1}{n} \sum
_{i=1}^{n}{\sigma_i \bigl(
\gamma(\c,X_i) - \gamma \bigl(\c^*,X_i\bigr)\bigr)}}
\\
&&\qquad \leq\sqrt{\frac{\pi}{2}}\mathbb{E}_{X,g} \sup
_{\llVert   \c- \c^* \rrVert  ^2 \leq\delta/16 M^2}{\frac{1}{n} \sum_{i=1}^{n}{g_i
\bigl(\gamma(\c,X_i) - \gamma\bigl(\c^*,X_i\bigr)
\bigr)}}.
\end{eqnarray*}
To derive bounds on the Gaussian complexity defined above, the
following comparison result between Gaussian processes is needed.

\begin{teo}[(Slepian's lemma)]\label{Slepian}
Let $X_t$ and $Z_t$, $t$ in $\mathcal{V}$, be some centered real
Gaussian processes. Assume that
\[
\forall s,t \in\mathcal{V}, \qquad\Var(Z_s - Z_t) \leq
\Var(X_s - X_t),
\]
then
\[
\mathbb{E} \sup_{t \in\mathcal{V}} Z_t \leq2 \mathbb{E} \sup
_{t
\in\mathcal{V}} X_t.
\]
\end{teo}

A proof of Theorem~\ref{Slepian} can be found in Theorem~3.14 of \cite
{Massart03}. For a fixed sample $X_1, \ldots, X_n$, define the
Gaussian process $Z_{\c}$ by
\[
Z_{\c} = \sum_{i=1}^{n}{g_i
\bigl(\gamma(\c,X_i) - \gamma\bigl(\c^*,X_i\bigr)
\bigr)},
\]
over the set $\mathcal{V}(\delta)=\mathcal{B}  (\c^*,\frac
{\sqrt{\delta}}{4M}  )$, where $\c^*$ is a fixed optimal
codebook. For $i = 1, \ldots, n$, $\c$, $\c' \in \mathcal
{V}(\delta)$, we have
\begin{eqnarray*}
\bigl( \gamma(\c,X_i) - \gamma\bigl(\c',X_i
\bigr) \bigr)^2 & \leq&\sup_{j=1, \ldots, k}{ \bigl( \llVert
X_i - c_j \rrVert ^2 - \bigl\llVert
X_i - c'_j \bigr\rrVert ^2
\bigr)^2}
\\
& \leq&\sup_{j=1, \ldots, k}{ \bigl( -2 \bigl\langle c_j -
c'_j, X_i \bigr\rangle+ \llVert
c_j\rrVert ^2 - \bigl\llVert c'_j
\bigr\rrVert ^2 \bigr)^2}
\\
& \leq&\sup_{j=1, \ldots, k}{ \bigl( 8 \bigl\langle c_j -
c'_j, X_i \bigr\rangle^2 + 2
\bigl(\llVert c_j\rrVert ^2 - \bigl\llVert
c'_j\bigr\rrVert ^2\bigr)^2
\bigr)}.
\end{eqnarray*}
Define now the Gaussian process $X_{\c}$ by
\[
X_{\c} = 2 \sqrt{2} \sum_{i=1}^{n}{
\sum_{j=1}^{k}{ \bigl\langle
c_j - c_j^*,X_i \bigr\rangle
\xi_{i,j}}} + \sqrt{2n} {\sum_{j=1}^{k}{
\bigl(\llVert c_j \rrVert ^2 - \bigl\llVert
c_j^*\bigr\rrVert ^2\bigr) \xi'_j}},
\]
where the $\xi$'s and $\xi'$'s are independent standard Gaussian
variables. It is straightforward that $\Var(Z_{\c} - Z_{\c'}) \leq
\Var(X_{\c} - X_{\c'})$. Therefore, applying Theorem~\ref{Slepian}
leads to
%
\begin{eqnarray} \label{decomposition}
\mathbb{E}_{g}{\sup_{\c\in\mathcal{V}(\delta)}{ Z_{\c} }} & \leq& 2 \mathbb{E}_{\xi}{\sup_{\c\in\mathcal{V}(\delta)}
{X_{\c} }}
\nonumber
\\
& \leq&4 \sqrt{2} \mathbb{E}_{\xi}{\sup_{\c\in\mathcal
{V}(\delta)} {\sum
_{i=1}^{n}{\sum
_{j=1}^{k}{ \bigl\langle c_j -
c^*_j,X_i \bigr\rangle\xi_{i,j}}}}}
\\
&&{}+ 2 \sqrt{2n} \mathbb{E}_{\xi'}{\sup_{\c\in\mathcal{V}(\delta
)}{\sum
_{j=1}^{k}{ \bigl(\llVert c_j
\rrVert ^2 - \bigl\llVert c^*_j \bigr\rrVert
^2\bigr) \xi'_{j}}}}.\nonumber
\end{eqnarray}
Using almost the same technique as in the proof of Theorem~2.1 in \cite
{Biau08}, the first term of the right-hand side of (\ref
{decomposition}) can be bounded as follows:
\begin{eqnarray*}
&& \mathbb{E}_{\xi}\sup_{\c\in\mathcal{V}(\delta)}{\sum
_{i=1}^{n}{\sum_{j=1}^{k}{
\bigl\langle c_j - c^*_j,X_i \bigr\rangle
\xi_{i,j}}}}
\\
&&\qquad  = \mathbb{E}_{\xi}{\sup_{\c\in\mathcal{V}(\delta
)}{
\sum_{j=1}^{k}{ \Biggl\langle
c_j - c^*_j, \Biggl( \sum_{i=1}^{n}{
\xi_{i,j} X_i } \Biggr) \Biggr\rangle}}}
\\
&&\qquad  \leq \mathbb{E}_{\xi}{\sup_{\c\in\mathcal{V}(\delta)}{\bigl\llVert \c-
\c^* \bigr\rrVert \sqrt{ \sum_{j=1}^{k}
{\Biggl\llVert \sum_{i=1}^{n}{
\xi_{i,j} X_i} \Biggr\rrVert ^2 }}}}
\\
&&\qquad  \leq\frac{\sqrt{\delta}}{4M} \sqrt{\sum_{j=1}^{k}
{\mathbb {E}_{\xi} \Biggl\llVert \sum_{i=1}^{n}{
\xi_{i,j} X_i} \Biggr\rrVert ^2 }}
\\
&&\qquad  \leq\frac{\sqrt{ k \delta}}{4M} \sqrt{\sum_{i=1}^{n}
{\llVert X_i \rrVert ^2}}.
\end{eqnarray*}
Then, applying Jensen's inequality ensures that
\[
\mathbb{E}_{X}\sqrt{\sum_{i=1}^{n}
{\llVert X_i \rrVert ^2}} \leq \sqrt{n} M.
\]
Similarly, the second term of the right-hand side of (\ref
{decomposition}) can be bounded from above by
\begin{eqnarray*}
&& \mathbb{E}_{\xi'}{\sup_{\c\in\mathcal{V}(\delta)}{\sum
_{j=1}^{k}{ \bigl(\llVert c_j \rrVert
^2 - \bigl\llVert c^*_j \bigr\rrVert ^2
\bigr) \xi'_{j}}}}
\\
&&\qquad  \leq\mathbb{E}_{\xi'}{\sup
_{\c\in\mathcal{V}(\delta)}{\sqrt {\sum_{j=1}^{k}
\bigl( \llVert c_j \rrVert ^2 - \bigl\llVert
c^*_j\bigr\rrVert ^2 \bigr)^2}\sqrt
{\sum_{j=1}^{k}{
\xi'_{j}} ^2}}}
\\
&&\qquad  \leq\frac{\sqrt{k \delta}}{2}.
\end{eqnarray*}
Combining these two bounds ensures that, for a fixed $\c^*$,
\[
\mathbb{E}_{X,g} \sup_{\llVert   \c- \c^*\rrVert  ^2 \leq\delta/ 16M^2}{Z_{\c
}}
\leq2 \sqrt{2kn} \sqrt{\delta},
\]
which leads to
%
\begin{equation}
\label{complexitecodebookfixe} \mathbb{E}Y_{\c^*} \leq\frac{4 \sqrt{k \pi\delta}}{\sqrt{n}}.
\end{equation}
Combining (\ref{deviationovercodebooks}) and (\ref
{complexitecodebookfixe}) into (\ref{decompositionwrtM}) gives the result.

\subsection{Proof of Proposition \protect\ref{minimax}}\label
{Proofofminimax}
Throughout this subsection, $\mathcal{H} = \mathbb{R}^d$, and for a
codebook $\c$, let $Q$ denote the associated nearest neighbor
quantizer. In the general case, such an association depends on how the
boundaries are allocated. However, since the distributions involved in
the minimax result have densities, how boundaries are allocated will
not matter.

Let $k \geq3$ be an integer. For convenience, $k$ is assumed to be
divisible by $3$. Let $m=2k/3 $. Let $z_1, \ldots, z_m$ denote a $6
\Delta$-net in $\mathcal{B}(0,M-\rho)$, where $\Delta>0$, and $w_1,
\ldots, w_m$ a sequence of vectors such that $\llVert  w_i\rrVert   = \Delta$.
Finally, denote by $U_i$ the ball $\mathcal{B}(z_i,\rho)$ and by
$U'_i$ the ball $\mathcal{B}(z_i + w_i,\rho)$. Slightly anticipating,
define $\rho= \frac{\Delta}{16}$.

To get the largest $\Delta$ such that for all $i=1,\ldots,m$, $U_i$
and $U'_i$ are included in $\mathcal{B}(0,M)$, it suffices to get the
largest $\Delta$ such that there exists a $6 \Delta$-net which can be packed in $\mathcal
{B}(0,M-\Delta/16)$. Since the cardinal of a maximal $6 \Delta$-net
is larger than the smallest number of balls of radius $6\Delta$ which
together cover $\mathcal{B}(0,M-\Delta/16)$, a sufficient condition
on $\Delta$ to guarantee that a $6\Delta$-net can be found is given by
\[
m \leq \biggl(\frac{M-\Delta/16}{6\Delta} \biggr)^d.
\]
Since $\Delta\leq M$, $\Delta$ can be chosen as
\[
\Delta= \frac{5M}{32 m^{1/d}}.
\]
For such a $\Delta$, $\rho$ takes the value $\rho= \frac{\Delta
}{16} = \frac{5M}{512 m^{1/d}}$. Therefore, it only depends on $k$,
$d$ and $M$.

Let $z=(z_i)_{i=1,\ldots,m}$ and $w=(w_i)_{i=1,\ldots,m}$ be
sequences as described above, such that, for $i=1,\ldots,m$, $U_i$ and
$U'_i$ are included in $\mathcal{B}(0,M)$.
For a fixed $\sigma\in \{-1,+1  \}^{m}$ such that $\sum_{i=1}^{m}{\sigma_i}=0$, let $P_\sigma$ be defined as
\[
\cases{ \displaystyle P_{\sigma}(U_i) = \frac{1+\sigma_i\delta}{2m},
\vspace*{3pt}\cr
\displaystyle P_{\sigma}\bigl(U'_i\bigr) =
\frac{1+\sigma_i\delta}{2m},
\vspace*{3pt}\cr
\displaystyle P_{\sigma} \underset{U_i} \sim\,
\bigl( \rho- \llVert x - z_i \rrVert \bigr)\mathbh{1}_{\llVert   x - z_i \rrVert   \leq\rho}
\,d\lambda(x),
\vspace*{3pt}\cr
\displaystyle P_{\sigma}\underset{U'_i} \sim\,\bigl( \rho- \llVert x - z_i - w_i\rrVert \bigr)
\mathbh{1}_{\llVert   x - z_i - w_i \rrVert   \leq\rho} \,d\lambda(x),}
\]
where $\lambda$ denotes the Lebesgue measure and $\delta\leq1/3$.
These distributions have been designed to have continuous cone-shaped
densities, as in Theorem~4 of \cite{Antos05}.

Similarly, let $Q_\sigma$ denote the quantizer defined by $Q_\sigma
(U_i) = Q_\sigma(U'_i) = z_i + \omega_i/2$ if $\sigma_i=-1$,
$Q_\sigma(U_i) = z_i$ and $Q_\sigma(U'_i) = z_i + \omega_i$ if
$\sigma_i = +1$. At last, for $\tau$ in $\{-1,+1\}^{m/2}$,
$\sigma(\tau)$ is defined as the sequence in $\{-1,+1\}^{m}$ such that
\[
\cases{\displaystyle \sigma_i(\tau)=\tau_i,
\vspace*{3pt}\cr
\displaystyle\sigma_{i+ m/2}(\tau)=-\sigma_i(\tau),}
\]
for $i=1,\ldots,\frac{m}{2}$, and the set of corresponding $Q_{\sigma
(\tau)}$'s is denoted by $\mathcal{Q}$.

Given a quantizer $Q$, let $R(Q,P_\sigma)$ and $\ell(Q,P_{\sigma})$
denote, respectively, the distortion and loss of $Q$ in the case where
the source distribution is $P_\sigma$. Proposition~\ref
{definitionDelta} below shows that only quantizers in $\mathcal{Q}$
may be considered in order to derive lower bounds on $R$.

%
\begin{prop}\label{definitionDelta}
Let $\sigma$ and $\sigma'$ be in $ \{-1,+1  \}^{{m}}$ such
that $\sum_{i=1}^{m}{\sigma_i} =\sum_{i=1}^{m}{\sigma'_i}=0$, and
let $\rho(\sigma,\sigma')$ denote the distance $\sum_{i=1}^{{m}}{\llvert  \sigma_i - \sigma'_i\rrvert  }$. Then
%
\begin{equation}
\label{distortionsigma} R(Q_{\sigma'}, P_\sigma) = R(Q_\sigma,
P_\sigma) + \frac{\Delta^2
\delta}{8m} \rho\bigl(\sigma,\sigma'
\bigr).
\end{equation}

Furthermore, for every nearest neighbor quantizer $Q$, there exists
$\sigma$ and $\tau$ such that
\[
\forall P_{\sigma(\tau')} \qquad R(Q,P_{\sigma(\tau')}) \geq
R(Q_\sigma,P_{\sigma(\tau')}) \geq\tfrac{1}{2}R(Q_{\sigma(\tau
)},P_{\sigma(\tau')}).
\]
At last, if $Q \neq Q_\sigma$, then the first inequality is strict,
for every $P_{\sigma(\tau')}$.
\end{prop}

The proof of Proposition~\ref{definitionDelta} follows the proof of
step~3 of Theorem~1 in \cite{Bartlett98}, and can be found in Section~B.1 (supplementary material \cite{supple}).

Since, for $\sigma\neq\sigma'$, $R(Q_\sigma',P_{\sigma}) >
R(Q_{\sigma},P_{\sigma})$, Proposition~\ref{definitionDelta} ensures
that the $P_{\sigma(\tau)}$'s have unique optimal codebooks, up to
relabeling. According to Proposition~\ref{definitionDelta}, the
minimax lower-bound over empirically designed quantizer may be reduced
to a lower-bound on empirically designed $\tau$'s, that is,
%
\begin{eqnarray}\label{minimaxreduction}
&& \inf_{\hat{Q}_n} \sup_{\tau\in\{-1,+1\}^{m/2}}{
\mathbb {E} \ell(\hat{Q}_n,P_{\sigma(\tau)})}\nonumber
\\
&&\qquad  \geq
\frac{1}{2}\inf_{\hat{\tau}} \sup_{\tau\in\{-1,+1\}^{m/2}}{
\mathbb{E} \ell(Q_{\sigma(\hat{\tau})},P_{\sigma(\tau)})}
\\
&&\qquad  \geq\frac{\Delta^2 \delta}{8m} \inf_{\hat{\tau}} \sup_{\tau
\in\{-1,+1\}^{m/2}}{
\mathbb{E} \rho(\hat{\tau},\tau)},\nonumber
\end{eqnarray}
where the inequality $\rho(\sigma(\tau),\sigma(\tau')) = 2 \rho
(\tau,\tau')$ has been used in the last inequality.

Let us define, for two distributions $P$ and $Q$ with densities $f$ and
$g$, the Hellinger distance
\[
H^2(P,Q)= \int_{\mathbb{R}^d} ( \sqrt{f} - \sqrt{g}
)^2(x) \,d\lambda(x).
\]
To apply Assouad's lemma to $\mathcal{Q}$, the following lemma is needed.

%
\begin{lem}\label{Assouadtechnique}
Let $\tau$ and $\tau'$ denote two sequences in $\{-1,+1\}^{m/2}$ such that $\rho(\tau,\tau') =2$. Then
\[
H^2\bigl(P^{\otimes n}_{\sigma(\tau)}, P^{\otimes n}_{\sigma(\tau')}
\bigr) \leq\frac{4n\delta^2}{m}:= \alpha,
\]
where $P^{\otimes n}$ denotes the product law of an $n$-sample drawn
from $P$.
\end{lem}

The proof of Lemma~\ref{Assouadtechnique} is given in Section~B.2 (supplementary material \cite{supple}).
Equipped with Lemma~\ref{Assouadtechnique}, a direct application
of Assouad's lemma as in Theorem~2.12 of \cite{Tsybakov09} yields,
provided that $\alpha\leq2$,
\[
\inf_{\hat{\tau}} \sup_{\tau\in\{-1,+1\}^{m/2}}{\mathbb {E} \rho(
\hat{\tau},\tau)} \geq\frac{m}{4} \bigl( 1 - \sqrt {\alpha ( 1 - \alpha/4
)} \bigr).
\]

Taking $\delta= \frac{\sqrt{m}}{2\sqrt{n}}$ ensures that $\alpha
\leq2$. For this value of $\delta$, it easily follows from~(\ref
{minimaxreduction}) that
\[
\sup_{\tau\in\{-1,+1\}^{m/2}}{\mathbb{E} \ell( \hat {Q}_n,P_{\sigma(\tau)}
)} \geq c_0 M^2 \sqrt{\frac{k^{1-4/d}}{n}},
\]
for any empirically designed quantizer $\hat{Q}_n$, where $c_0$ is an
explicit constant.

Finally, note that, for every $\delta\leq\frac{1}{3}$ and $\sigma$,
$P_\sigma$ satisfies a margin condition as in~(\ref
{conditionmargeminimax}), and is $\varepsilon$-separated, with
\[
\varepsilon= \frac{\Delta^2 \delta}{2m}.
\]
This completes the proof of Proposition~\ref{minimax}.

\subsection{Proof of Proposition \protect\ref{Gaussianboundary}}

As mentioned below in Proposition~\ref{Gaussianboundary}, the inequality
\[
\frac{\theta_{\min}}{\theta_{\max}} \geq\frac{2048 k \sigma
^2}{(1-\varepsilon) \tilde{B}^2(1 - e^{-\tilde{B}^2/{2048\sigma^2}})}
\]
ensures that, for every $j$ in $ \{ 1,\ldots, k  \}$, there
exists $i$ in $ \{ 1,\ldots, k  \}$ such that $\llVert
c_i^* - m_j \rrVert   \leq\tilde{B}/16$. To be more precise, let $
\mathbf{m}$ denote the vector of means $(m_1, \ldots, m_k)$, then
\begin{eqnarray*}
R(\mathbf{m}) & \leq&\sum_{i=1}^{k}{
\frac{\theta_i}{2 \pi\sigma^2
N_i} \int_{V_i(\mathbf{m})}{\llVert x - m_i
\rrVert ^2 e^{-{\llVert   x - m_i \rrVert
^2}/(2 \sigma^2)}\,d \lambda(x)}}
\\
& \leq&\frac{p_{\max}}{2(1-\varepsilon) \pi\sigma^2} \sum_{i=1}^{k}{
\int_{\mathbb{R}^2}{\llVert x - m_i\rrVert ^2
e^{-{\llVert   x - m_i \rrVert
^2}/(2 \sigma^2)}\,d \lambda(x)}}
\\
& \leq&\frac{2kp_{\max} \sigma^2}{1-\varepsilon}.
\end{eqnarray*}
Assume that there exists $i$ in $ \{1, \ldots, k  \}$ such
that, for all $j$, $\llVert   c_j^* - m_i \rrVert   \geq\tilde{B}/16$. Then
\begin{eqnarray*}
R(\c) & \geq&\frac{\theta_i}{2 \pi\sigma^2}\int_{\mathcal
{B}(m_i,\tilde{B}/32)}{\frac{\tilde{B}^2}{1024}
e^{-{\llVert  x -
m_i\rrVert  ^2}/(2 \sigma^2)}\,d \lambda(x)}
\\
& \geq&\frac{\tilde{B}^2 \theta_{\min}}{2048 \pi\sigma^2} \int_{\mathcal{B}(m_i,\tilde{B}/32)}{e^{-{\llVert  x - m_i\rrVert  ^2}/(2 \sigma
^2)}\,d
\lambda(x)}
\\
& >& \frac{\tilde{B}^2 \theta_{\min}}{1024} \bigl( 1 - e^{-{\tilde{B}^2}/(2048 \sigma^2)} \bigr)
\\
& >& R(\mathbf{m}).
\end{eqnarray*}
Hence, the contradiction. Up to relabeling, it is now assumed that for
$i=1,\ldots,k$, $\llVert  m_i - c_i^*\rrVert   \leq\tilde{B}/16$. Take $y$ in
$N_{\c^*}(x)$, for some $\c^*$ in $\mathcal{M}$ and for $x \leq
\frac{\tilde{B}}{8}$, then, for every $i$ in $ \{ 1, \ldots, k
 \}$,
\[
\llVert y - m_i \rrVert \geq\frac{\tilde{B}}{4},
\]
which leads to
\[
\sum_{i=1}^{k}{\frac{\theta_i}{2 \pi\sigma^2 N_i}{\llVert
y - m_i\rrVert ^2 e^{-{\llVert   y - m_i \rrVert  ^2}/(2 \sigma^2)}}} \leq
\frac{k \theta
_{\max}}{(1-\varepsilon) 2 \pi\sigma^2} e^{-{\tilde{B}^2}/(32
\sigma^2)}.
\]
Since the Lebesgue measure of $N_{\c^*}(x)$ is smaller than $4k \pi M
x$, it follows that
\[
P\bigl(N_{\c^*}(x)\bigr) \leq\frac{2 k^2 M \theta_{\max}}{(1-\varepsilon)
\sigma^2} e^{-{\tilde{B}^2}/(32 \sigma^2)}x.
\]
On the other hand, $\llVert  m_i - c_i^*\rrVert   \leq\tilde{B}/16$ yields
\[
\mathcal{B}(m_i,3 \tilde{B} /8) \subset V_i\bigl(\c^*
\bigr).
\]
Therefore,
\begin{eqnarray*}
P\bigl(V_i\bigl(\c^*\bigr)\bigr) & \geq&\frac{\theta_i}{2 \pi\sigma^2 N_i} \int
_{\mathcal{B}(m_i,3 \tilde{B} /8)}{e^{-{\llVert  x - m_i \rrVert  ^2}/(2
\sigma^2)}\,d \lambda(x)}
\\
& \geq&\theta_i \bigl(1 - e^{-{9 \tilde{B}^2}/(128 \sigma^2)} \bigr),
\end{eqnarray*}
hence $p_{\min} \geq\theta_{\min} (1 - e^{-{9 \tilde
{B}^2}/(128 \sigma^2)}  )$. Consequently, provided that
\[
\frac{\theta_{\min}}{\theta_{\max}} \geq\frac{2048 k^2 M^3}{(1 -
\varepsilon)7\sigma^2 \tilde{B}(e^{\tilde{B}^2/{32\sigma^2}}-1)},
\]
direct calculation shows that
\[
P\bigl(N_{\c^*}(x)\bigr) \leq\frac{B p_{\min}}{128 M^2}x.
\]
This ensures that $P$ satisfies (\ref{majorationkappa}). According to
(ii) in Proposition~\ref{lienconditionmarge}, $\mathcal{M}$ is finite.

\section*{Acknowledgments} The author would like to thank three
referees for valuable
comments and suggestions.

\begin{supplement}[id=suppA]
\sname{Appendix}
\stitle{Remaining proofs}
\slink[doi]{10.1214/14-AOS1293SUPP} 
\sdatatype{.pdf}
\sfilename{aos1293\_supp.pdf}
\sdescription{Due to space constraints, we relegate technical details
of the remaining proofs to the supplement \cite{supple}.}
\end{supplement}

%

\printaddresses
\end{document}